\documentclass[a4paper,11pt]{article}

\usepackage{amsmath}
\usepackage{amssymb}
\usepackage{amsthm}
\usepackage{graphicx}
\usepackage{amscd}
\usepackage{epic, eepic}
\usepackage{url}
\usepackage{color}
\usepackage[utf8]{inputenc} 
\usepackage{comment}
\usepackage{enumerate}   
\usepackage{graphicx}
\usepackage{epstopdf}
\usepackage{enumitem}
\usepackage{tikz}
\usepackage[top=24mm, bottom=25mm, left=25mm, right=25mm]{geometry}
\usepackage[bookmarks=false,hidelinks]{hyperref}
\usepackage{mathrsfs}
\usetikzlibrary{arrows}
\definecolor{wrwrwr}{rgb}{0.3803921568627451,0.3803921568627451,0.3803921568627451}
\newcommand{\abs}[1]{\left\vert{#1}\right\vert}

\newtheorem{theorem}{Theorem}
\newtheorem{claim}[theorem]{Claim}

\newtheorem{lemma}[theorem]{Lemma}
\newtheorem{observation}[theorem]{Observation}

\newcommand{\ex}{{\rm  ex}}

\usepackage{authblk}

\title{Exact results for generalized extremal problems forbidding an even cycle}
\author[1]{Ervin Győri}
\author[1,3]{Zhen He}
\author[1,3]{Zequn Lv}
\author[1]{Nika Salia}
\author[1]{Casey Tompkins}
\author[1,4]{Kitti Varga}
\author[1,2]{Xiutao Zhu}
\date{}

\affil[1]{Alfr\'ed R\'enyi Institute of Mathematics, Hungarian Academy of Sciences. }
\affil[2]{Department of Mathematics, Nanjing University.}
\affil[3]{Department of Mathematical Sciences, Tsinghua University.}
\affil[4]{Department of Computer Science and Information Theory, Budapest University of Technology and Economics.}

\begin{document}

\maketitle
\begin{abstract}
We determine the maximum number of copies of $K_{s,s}$ in a $C_{2s+2}$-free $n$-vertex graph for all  integers $s \ge 2$ and sufficiently large $n$.  Moreover, for $s\in\{2,3\}$ and any integer $n$ we obtain the maximum number  of cycles of length $2s$ in an  $n$-vertex  $C_{2s+2}$-free bipartite graph.
\end{abstract}

\section{Introduction}
 \textbf{Notation.} For a graph $G$, the vertex and edge set of $G$ are denoted by $V(G)$ and $E(G)$, respectively.
 Furthermore, we let $v(G)=\abs{V(G)}$ and $e(G)=\abs{E(G)}$. The maximum degree of a vertex in $G$ is denoted by $\Delta(G)$.
For $S \subseteq V(G)$, we denote by $G[S]$ the subgraph of $G$ induced on the vertex set $S$.  
For a vertex $v$, we denote by $G-v$ the graph $G[V(G)\setminus \{v\}]$.
Moreover, for a subgraph $H$ of $G$ we write $G-H$ for $G[V(G)\setminus V(H)]$.
For $v\in V(G)$, we denote by $N(v)$ the neighborhood of $v$ and we write $N_H(v)$ for $N(v)\cap V(H)$. For a vertex $v$, the set $N(v) \cup \{v\}$ is denoted by~$N[v]$. 
For a set $S$, we denote the set of $k$-element subsets of $S$ by $\binom{S}{k}$.
By a copy of~$H$ in a graph $G$ we mean a subgraph of $G$ isomorphic to $H$. 
For graphs $G$ and $H$, we denote by $G+H$ the join of $G$ and $H$, that is the graph obtained by connecting each pair of vertices between a vertex disjoint copy of $G$ and $H$. 

The path, cycle and complete graph on $k$ vertices are denoted by $P_k$, $C_k$ and $K_k$, respectively. Let $S_k$ be a $k$-vertex complete bipartite graph with a color class of size $1$.
We sometimes refer to a $k$-vertex path by a sequence of vertices $v_1 v_2 \dots v_k$, where $v_i v_{i+1}$ is an edge for all $1 \le i \le k-1$, and we sometimes refer to a $k$-vertex cycle by a sequence of the form $v_1 v_2 \dots v_k v_1$ where $v_1 v_k$ and  $v_i v_{i+1}$ are edges for all $1 \le i \le k-1$.
We denote by $K_{a,b}$ the complete bipartite graph with parts of size $a$ and $b$, and we denote by $K_{a,b}^+$ the join of a clique of size $a$ and an independent set of size $b$. 
A star is a tree with at most one vertex of degree at least two; this vertex is referred to as the central vertex of the star. A forest in which every connected component is a star is called a star forest. A double star is a tree consisting of two adjacent vertices $u$ and $v$ such that every vertex in $(N(u)\setminus \{v\}) \cup (N(v)\setminus \{u\})$ is a leaf; such vertices $u$ and $v$ are referred to as central vertices of the double star.
For any two graphs $H$ and $G$, we denote the number of copies of $H$ in $G$ by $H(G)$.
For given graphs $F$ and $G$, we say that $G$ is $F$-free if it does not contain $F$ as a  subgraph (not necessarily induced). 
The maximum of $H(G)$ across $n$-vertex $F$-free graphs $G$ is denoted by $\ex(n,H,F)$. 
The maximum of $H(G)$ across bipartite $n$-vertex $F$-free graphs $G$ is denoted by $\ex_{bip}(n,H,F)$. 
\vspace{.4em}

\noindent \textbf{Background.}  Generalized extremal problems have a long history dating back to a result of Zykov~\cite{zykov1949some} (and later independently Erd\H{o}s~\cite{erdos1962number}), who determined for all $s$ and $t$ the value of $\ex(n,K_s,K_t)$, extending the classical theorem of Tur\'an~\cite{turan1941external}. 
Following this initial result, there has been extensive work determining the value of $\ex(n,H,F)$ for various pairs of graphs $H$ and $F$.  
An important early result in this direction is a result of Gy\H{o}ri, Pach and Simonovits~\cite{gyori1991maximal} which determined the value of $\ex(n,H,K_3)$ for all bipartite graphs $H$ containing a matching on all but at most one of its vertices.
The generalized extremal function $\ex(n,H,F)$ for arbitrary pairs of graphs was introduced by Alon and Shikhelman~\cite{alon2016many}. 
In addition to the function $\ex(n,H,F)$, we will consider the analogous function under the assumption that the ground graph is bipartite, which we denote by $\ex_{bip}(n,H,F)$.

In the present paper we will be interested in the setting when the forbidden graph is an even cycle.
Determining $\ex(n,C_5,C_3)$ was a long standing open problem of Erd\H{o}s~\cite{erdos1984some}.
A construction is obtained by taking a blow-up of $C_5$ consisting of almost equal classes. This problem was finally settled by Grzesik~\cite{grzesik2012maximum} and independently by Hatami, Hladk{\`y}, Kr{\'a}l, Norine and Razborov~\cite{hatami2013number}. This result has recently been extended Grzesik and Kielak~\cite{grzesik2022maximum} who determined $\ex(n,C_{2k+1},C_{2k-1})$ asymptotically for all $k>2$.  
The dual problem of determining $\ex(n,C_3,C_5)$ was introduced by Bollob\'as and Gy\H{o}ri~\cite{bollobas2008pentagons} who proved an upper bound, and this upper bound was subsequently improved in~\cite{alon2016many},~\cite{ergemlidze2019note} and~\cite{ergemlidze2018triangles}. Bounds on $\ex(n,C_3,C_{2k+1})$ were obtained by Gy\H{o}ri and Li~\cite{gyHori2012maximum}.
An exact result for $\ex(n,C_5,C_7)$ was obtained by G\'orski and Grzesik~\cite{Gorski}.
The order of magnitude of the $\ex(n,C_k,C_\ell)$ for general $k$ and $\ell$ was determined by Gishboliner and Shapira~\cite{gishboliner2020generalized} and independently for even $k$ and $\ell$ by Gerbner, Gy\H{o}ri, Methuku and Vizer~\cite{gerbner2020generalized}. 

In~\cite{gerbner2020generalized} an asymptotic result was obtained for $\ex(n,C_4,C_6)$ and $\ex_{bip}(n,C_6,C_8)$.  We determine the exact value of $\ex_{bip}(n,C_4,C_6)$, $\ex_{bip}(n,C_6,C_8)$ for all $n$, and the exact value of $\ex(n,C_4,C_6)$ for sufficiently large $n$ in the following three theorems. Moreover, we determine the structure of the extremal graphs in each case.   

\begin{theorem}\label{thm:bipC4C6}
For all positive integers $n$, we have
\[
\ex_{bip}(n,C_4,C_6) = \binom{n-2}{2},
\]
and equality holds only for $K_{2,n-2}$.
\end{theorem}

\begin{theorem}\label{thm:bipC6C8}
For all positive integers $n$, we have
\[
\ex_{bip}(n,C_6,C_8) = 6 \binom{n-3}{3},
\]
and equality holds only for $K_{3,n-3}$.
\end{theorem}



\begin{theorem}\label{thm:C_4,C_6-free}
For $n>3\binom{31}{4}$, we have
\[
\ex(n,C_4,C_6)=\binom{n-2}{2}+2,
\]
and equality holds only for a graph obtained from $K_{2,n-2}$ by adding a single edge in the both independent sets.
\end{theorem}

Finally, we turn our attention to exact results about the generalized extremal numbers of complete bipartite graphs $K_{s,s}$ in graphs without a copy of $C_{2s+2}$.  We prove the following exact result.

\begin{theorem} \label{thm:kss}
For integers $n$ and $s$ such that $s\geq 3$ and $n\geq 3s+1+\binom{2s+1}{s}\frac{s+1}{2}$, we have
\[
\ex(n,K_{s,s},C_{2s+2})=\binom{n-s}{s}.
\]
Equality holds only for graphs containing $K_{s,n-s}$ and contained in the graph $K_s+H$, where $H$ is an $(n-s)$-vertex graph with exactly one edge.
\end{theorem}

The proofs of Theorems~\ref{thm:bipC4C6}~and~\ref{thm:bipC6C8} are given in Section~\ref{sec:bip}, and the proofs of Theorems~\ref{thm:C_4,C_6-free}~and~\ref{thm:kss} are given in Section~\ref{sec:rest}.

\section{Proofs of theorems about bipartite graphs} \label{sec:bip}

\begin{proof}[Proof of Theorem \ref{thm:bipC4C6}]
Consider the graph $K_{2,n-2}$. It contains no cycle of length six and $ \binom{n-2}{2}$ cycles of length four. Hence we have $\ex_{bip}(n,C_{4},C_{6})\geq \binom{n-2}{2}$. 

Let $G$ be an $n$-vertex $C_6$-free bipartite graph. 
Let $C$ be a cycle of length four in $G$ such that  $C=v_1 v_2 v_3 v_4 v_1$.
Then since $G$ is a bipartite $C_6$-free graph, either there is no cycle of length four containing the vertices $v_1,v_3$  distinct from $C$ or there is no cycle of length four containing the vertices $v_2,v_4$  distinct from $C$. Hence each $C_4$ in $G$ contains a pair of vertices from the same color class of $G$ which is unique for that cycle.  

We define an auxiliary graph $G'$ on the same set of vertices as $G$ by taking the edge $uv$ in $G'$ if and only if  there is exactly one $C_4$ in $G$ with opposite vertices $u$ and $v$.
Note that since $G$ is bipartite, the vertices $u$ and $v$ belong to the same color class of $G$. 
The number of copies of $C_4$ in $G$ is at most the number of edges in $G'$ since each copy of $C_4$ in $G$ contains a pair of opposite vertices which is unique for the cycle.  
If the color classes of $G$ have sizes $a$ and $n-a$ such that $2\leq a\leq \frac{n}{2}$, then the number of edges of $G'$ is at most $\binom{a}{2}+\binom{n-a}{2}$.
Note that if $a>2$, then $\binom{a}{2}+\binom{n-a}{2}<\binom{n-2}{2}$ when $n>6$.
Therefore if $n>6$ and $G$ contains at least $\binom{n-2}{2}$ copies of $C_4$, then it is a bipartite graph with color classes of sizes $n-2$ and $2$, and it follows that the number of copies of $C_4$ in $G$ is at most $\binom{n-2}{2}$ with equality only when $G$ is isomorphic to $K_{2,n-2}$.
For $n\leq 6$ it is straightforward to verify that the theorem holds.
\end{proof}

\begin{proof}[Proof of Theorem \ref{thm:bipC6C8}]
The graph $K_{3,n-3}$ contains no cycle of length eight and contains $ 6\binom{n-3}{3}$ cycles of length six. Hence we have $\ex_{bip}(n,C_{6},C_{8})\geq 6\binom{n-3}{3}$. 
 
Let $G$ be an $n$-vertex bipartite graph with color classes $V_1$ and $V_2$, containing no cycle of length eight and containing the maximum number of cycles of length six. 
Note that we may assume that $n\geq8$ since the theorem is simple to verify for $n<8$.
In what follows we are going to find a function $\phi$ from the set of copies of $C_6$ in $G$ to unordered monochromatic triples of the vertices of $G$. 
Then we are going to show that the inverse image of each triple has size at most six. This implies that 
\[
C_6(G)<6\binom{\abs{V_1}}{3}+6\binom{\abs{V_2}}{3}.
\]
If  $\abs{V_1},\abs{V_2}> 3$, then  we have $C_6(G)\leq 6\binom{\abs{V_1}}{3}+6\binom{\abs{V_2}}{3}< 6\binom{n-3}{3}$ since $n\geq 8$, a contradiction.
We may therefore assume $\abs{V_1}=3$, and it follows that $G$ is isomorphic to  $K_{3,n-3}$ since it maximizes the number of cycles of length six. Hence we are done if such a function $\phi$ exists.

We are going to define the function $\phi$ recursively. 
Let $C$ be a subgraph of $G$ be isomorphic to~$C_6$, such that $C=v_1 v_2 v_3 v_4 v_5 v_6 v_1$ with color classes $B=\{v_1,v_3,v_5\}$ and $R=\{v_2,v_4,v_6\}$, and assume that $\phi(C)$ is not yet defined for $C$.
If all of the cycles of length six containing $B$ as a color class contain $R$ as the other color class, then we take $\phi(C):=B$;
if all of the cycles of length six containing $R$ as a color class contain $B$ as the other color class, then we take $\phi(C):=R$. If both of the above occur, then we take $\phi(C):=B$. 
Otherwise, without loss of generality, we may assume there is a vertex~$v_2'$ incident to~$v_1$ and~$v_3$, and there is a vertex~$v_5'$ incident to~$v_4$ and~$v_6$. 

Consider a subgraph $G'$ of $G$ containing all cycles of length six with opposite vertices $v_1$ and $v_4$ such that each edge of $G'$ is in a cycle of length six with opposite vertices $v_1$ and $v_4$.
Let the neighborhood of $v_1$ in $G'$ excluding $v_4$ (if they are incident) be denoted by $U$, and let the neighborhood of $v_4$ in $G'$ excluding $v_1$ be denoted by $W$. Note that  $\abs{U},\abs{W}\geq 3$ since $\{v_2,v_2',v_6\}\subseteq U$ and $\{v_3,v_5,v_5'\}\subseteq W$.
Since $G'$ is $C_8$-free, it is straightforward to deduce that the structure of $G'$ is one of the following:
\begin{itemize}
    \item If $v_3v_6$ is an edge, then  $v_1v_4$ is also an edge since $G'$ is $C_8$-free and  every edge is in a copy of~$C_6$. 
    Then the induced graph $G'[U\cup W]$ is a double star with central vertices $v_3$ and $v_6$. 
    Note that every cycle of length six in $G'$ contains a vertex from $U\setminus \{v_6\}$ and a vertex from  $W\setminus \{v_3\}$.
    Even more, for each pair consisting of a vertex from $U\setminus \{v_6\}$ and a vertex from  $W\setminus \{v_3\}$, there are exactly two cycles of length six in $G'$. 
    Hence the number of cycles of length six in $G'$ is   $2(\abs{U}-1)(\abs{W}-1)$.
  
    \item  Otherwise  $v_3v_6$ is not an edge of $G'$, and  the  graph $G'[U\cup W]$ is a star forest. In particular, it contains at least two stars, one with central vertex $v_3$ and the other with central vertex $v_6$.
    Note that every~$C_6$ in~$G'$ contains $v_1$ and $v_4$ and at least two leaves of $G'[U\cup W]$. 
    Each pair of the leaves is contained in at most one cycle of length six in $G'$. 
    Hence the number of cycles of length six in $G'$ is  at most $\binom{\abs{U}+\abs{W}-2}{2}$.
\end{itemize}

Therefore we have the following bound on the number of cycles of length six in $G'$:
\[C_6(G')\leq \mbox{max}\left(\binom{\abs{U}+\abs{W}-2}{2} ,2(\abs{U}-1)(\abs{W}-1)\right).\]

Let us consider the set of triples
\[
X_{v_1,v_4}:=\left\{ \{v_1,w_1,w_2\}:\{w_1,w_2\}\subseteq W \right\}\cup \left\{\{v_4,u_1,u_2\}: \{u_1,u_2\}\subseteq U \right\}.
\]

Note that $\abs{X_{v_1,v_4}}=\binom{\abs{U}}{2}+\binom{\abs{W}}{2}$, hence $ 2\abs{X_{v_1,v_4}}>C_6(G')$.

\begin{observation}\label{Observation:cherry_structure_3-path}
There is no vertex in $V(G-G')$ incident to two vertices of $G'$ since $G$ is $C_8$-free. Hence the only vertices  adjacent to more than one vertex from $\{v_1\} \cup W$ are  $v_4$ and $v_6$.

There is no path of length three between any two vertices of $G'$ containing a vertex outside of $G'$ since $G$ is $C_8$-free and there is no vertex in $V(G-G')$ incident to two vertices of $G'$.

By the previous two claims, $v_6$ and $w\in W$ do not appear as opposite vertices in any copy of $C_6$ which is not contained in $G'$. 
Hence if $X_{v_1,v_4}\cap X_{u,v}\neq \emptyset$ for some opposite pair of vertices $u,v$,  then every six-cycle with opposite vertices $u$ and $v$ are six cycles of $G'$. Hence $G'$ doesn't contain any $6$-cycle  for which $\phi$  is already defined.
\end{observation}

Let  us denote the set of $6$-cycles in $G'$  by $\mathcal{C}$. We define $\phi$ on $\mathcal{C}$ in the following way: $\phi(\mathcal{C})\subseteq X_{v_1,v_4}$ and  $\abs{\phi^{-1}(x)\cap \mathcal{C}}\leq 2$ for every $x\in X_{v_1,v_4}$. 

We repeat the procedure iteratively for all cycles of length six in $G$ for which $\phi$ is not yet defined. 
By Observation \ref{Observation:cherry_structure_3-path}, for each triple $\{v_1,w_1,w_2\}\in X_{v_1,v_4}$ we have $\abs{\phi^{-1}(\{v_1,w_1,w_2\})}\leq 2<6$. 
\end{proof}

\section{Proofs of theorems about general graphs} \label{sec:rest}

\subsection{Maximizing cycles of length four}

First we prove an essential lemma which we will need in the proof of Theorem~\ref{thm:C_4,C_6-free}.

\begin{lemma}\label{lemma:c_4vertex}
Let $G$ be a $K_{2,s+1}$-free and $C_6$-free graph for some integer $s \ge 2$, and let $v\in V(G)$. 
Then the number of copies of $C_4$ incident to $v$ is at most
\[
\max \left(3\abs{N_{G}(v)}, \frac{(s-1)(s+2)}{2(s+1)}\abs{N_{G}(v)} \right).
\]
If $G$ is also $K_5$-free, then the number of copies of $C_4$ incident to $v$ is at most 
\[
\max \left(2\abs{N_{G}(v)},  \frac{(s-1)(s+2)}{2(s+1)}\abs{N_{G}(v)} \right).
\]
\end{lemma}

\begin{proof}
Let $G'$ be the subgraph of $G$ spanned by the edges which are in a copy of $C_4$ incident to $v$. 
By the linearity of the desired upper bound in Lemma~\ref{lemma:c_4vertex}, we may assume that $G'$ is $2$-connected. 
Let us denote the degree of $v$ in $G'$ by $x$.

First we consider the case when $v(G')=x+1$.
Observe that $G'-v$ is a connected graph since $G'$ is $2$-connected, and moreover we have that $G'[N_{G'}(v)]=G'-v$.  
Note that the number of paths of length two in $G'-v$ is equal to the number of copies of $C_4$ incident to $v$ in~$G$. 
The graph $G'-v$ is $P_5$-free since $G'$ is $C_6$-free. 
Even more, we have $\Delta(G'-v)\leq s$ since $G$ is $K_{2,s+1}$-free.
If $G'-v$ is a tree, then it is either a star or a double star, and therefore the number of copies of $C_4$ incident to $v$ in $G$ is at most $\frac{(s-1)x}{2}<\frac{(s-1)(s+2)}{2(s+1)}\abs{N_{G}(v)}$. 
If $G'-v$ contains a cycle of length four, then the number of copies of $C_4$ incident to $v$ is at most $12\leq 3x$. 
Note that if $G$ is $K_5$-free (and $G'-v$ contains a $C_4$), then the number of copies of $C_4$ incident to $v$ is at most $8\leq 2x$.
If $G'-v$ contains a triangle and no $4$-cycle, then $G'-v$ is a graph obtained from a triangle by adding $x-3$ pendant edges to one of its vertices.
In this case we have at most $\binom{x-1}{2}+2$ cycles of length~$4$ incident to $v$. 
Note that $x-1\leq s$, and a simple calculation shows that the required upper bound in Lemma~\ref{lemma:c_4vertex} holds.

From here we will assume that $v(G')>x+1$.
Hence $V(G')\setminus N_{G'}[v]$ is nonempty. 
Note that $G'[V(G')\setminus N_{G'}[v]]$ is an independent set, and every vertex from $V(G')\setminus N_{G'}[v]$ has at least two neighbors in $N_{G'}(v)$. 
Even more, if  $\abs{V(G')\setminus N_{G'}[v]}>1$, then for all vertices $u,u'\in V(G')\setminus N_{G'}[v]$ we have $N_{G'}(u)=N_{G'}(u')$ and $\abs{N_{G'}(u)}=2$ since $G'$ is a $2$-connected $C_6$-free graph.

If $G'[N_{G'}(v)]$ is connected, then $G'[N_{G'}(v)]$ is either a star or a triangle. 
It follows that the number of copies of $C_4$ incident to $v$ in $G'$ is at most 
\[
\max \left(2x, \frac{(s-1)(s+2)}{2(s+1)}x \right).
\]
Note that if $G'[N_{G'}(v)]$ is a star of $s+1$ vertices and  $G'[V(G')\setminus N_{G'}[v]]$ is an independent set of size $s-1$ incident to the same two vertices in $G'[N_{G'}(v)]$ inducing an edge, then there are exactly $\frac{(s-1)(s+2)}{2(s+1)}x$ copies of $C_4$ incident to $v$.

If $N_{G'}(v)$ is not connected, then  $G'[N_{G'}(v)]$ consists of isolated vertices and a star.
If  $G'[V(G')\setminus N_{G'}[v]]$ consists of just a vertex,  denoted by $u$, then in $G'$ all of the  $y$ isolated vertices are incident to $u$, therefore  $y<s$ and $x-y-1\leq s$. 
The star $S$ in $G'[N_{G'}(v)]$ is isomorphic to $S_{x-y}$. Since  every edge of $G'$ is in a copy of $C_4$ incident to $v$, the star $S$ is not an edge. The only vertex of $S$ incident to $u$ is the central vertex of the star. 
Then the number of copies of $C_4$ incident to $v$ is $\binom{y+1}{2}+\binom{x-y-1}{2}\leq \frac{(s-1)x}{2}< \max \left(2x,  \frac{(s-1)(s+2)}{2(s+1)}x \right)$ since $x-y-2\leq s-1$ and $y\leq s-1$. 
If $G'[V(G')\setminus N_{G'}[v]]$ contains at least two vertices, then each of the vertices have at least two neighbors in $N_{G'}(v)$. Hence $G'[N_{G'}(v)]$ contains  two connected components  since $G'$ is $2$-connected. Therefore 
$G'[N_{G'}(v)]$ consists of an isolated vertex $w$ and a star isomorphic to $S_{x-1}$ with central vertex~$w'$. Note that $S_{x-1}$ is not an edge and every vertex of $G'[V(G')\setminus N_{G'}[v]]$ is incident to $w$ and $w'$.   Then the number of copies of $C_4$ is at most $s-1+\binom{x-2}{2}\leq \frac{(s-1)x}{2}< \max \left(2x,  \frac{(s-1)(s+2)}{2(s+1)}x \right)$.
 \end{proof}

\begin{proof}[Proof of Theorem~\ref{thm:C_4,C_6-free}]
Let $G$ be an $n$-vertex graph with no cycle of length six. 
Let $K$ be a subgraph of $G$ isomorphic to $K_{2,s}$, where $s$ is the maximum integer for which $G$ contains a copy of $K_{2,s}$. 
Let the color classes of $K$ of size $2$ and $s$ be $A$ and $B$, respectively.

\begin{observation} \label{main}
For $s>2$, let $v$ and $u$ be distinct vertices of $V(G-K)$. Then the following properties hold. 

\begin{enumerate}
\item \label{2A} The vertex $v$ is not adjacent to both vertices in $A$ by the maximality of $s$.
\item \label{2B} Since $G$ is $C_6$-free, $v$ is not adjacent to more than one vertex from $B$. 
\item \label{2nbrs} By Items~\ref{2A} and~\ref{2B}, we have that $v$ and $u$ are adjacent to at most $2$ vertices of $K$. 
If~$v$ and~$u$ both have two neighbors in $V(K)$, then $N(v)\cap V(K)=N(u)\cap V(K)$ since $G$ is  $C_6$-free.
\item \label{L} For $s>3$, the graph $G[B]$ contains at most one edge since $G$ is $C_6$-free.
\end{enumerate}
\end{observation}

In the following proof we need to classify the copies of $C_4$ depending how they meet with~$K$. 
For simplicity, we introduce the following technical definitions.
If there exists a vertex~$v \in V(G-K)$ with two neighbors $x$ and $y$ in $V(K)$, then we call $xy$ a red edge. 
Note that there is at most one red edge in $G[V(K)]$ by Item~\ref{2nbrs} of Observation~\ref{main}. 
We refer to any edge in $G[A]$ or $G[B]$ as a blue edge.  
Note that there are at most two blue edges in $G[V(K)]$ for $s>3$ by Item~\ref{L} of Observation~\ref{main}. .

Let $C$ be a subgraph of $G$ isomorphic to $C_4$.
\begin{itemize}
    \item We say $C$ is of \emph{Type~$1$} if $\abs{V(K)\cap V(C)}=1$.
    \item We say $C$ is of \emph{Type~$2$} if $\abs{V(K)\cap V(C)}=2$.
    We will further divide the class of Type~$2$ copies of $C_4$ into two subclasses.
    Each copy of $C_4$ using a blue edge of $G[V(K)]$ will be referred to as a \emph{blue Type~$2$} $C_4$, and each copy of $C_4$ with two opposite vertices belonging to $V(K)$ (thus, inducing a red edge in $K$ as its diagonal) will be referred to as a \emph{red Type~$2$} $C_4$. Note that there is no $C_4$ of Type~$2$ which is neither blue nor red.
    
    \item We say $C$ is of \emph{Type~$3$} if $\abs{V(K)\cap V(C)}=3$. Note that such a copy of $C_4$ is incident to a blue edge, and it induces a red edge of $K$ as its diagonal.
    \item We say $C$ is of \emph{Type~$4$} if $\abs{V(K)\cap V(C)}=4$.
\end{itemize}

\begin{observation}\label{Obs:general_observation_red_blue}
It is important to note that if there is a blue Type~$2$  $C_4$, then $K$ contains exactly one blue edge since $G$ is $C_6$-free.
Thus there are at most $s-1$ blue Type~$2$ $C_4$'s since $G$ is $C_6$- and $K_{2,s+1}$-free.
Even more, if there is a blue edge and a red edge, then the blue edge is in the color class $A$. 
All Type~$3$ $C_4$'s share the same blue edge, hence they share the same three vertices of $K$. 
Therefore the number of Type~$3$ $C_4$'s is at most $s-1$. 
\end{observation}

In the following we show that for all $n\geq 31$ either we have $C_4(G)\leq \binom{n-2}{2}+2$, and equality holds if and only if $G$ is isomorphic to $K_{2,n-2}^+$ or $C_4(G)<\ex(n-1,C_4,C_6)+n-3$.
For a vertex~$v$ and for a copy~$K$ of $K_{2,s}$ with $v\in V(K)$, the collection of copies of $C_4$ incident to~$v$ in $G$ can be partitioned into five sets: Type~$4$, Type~$3$, red Type~$2$, blue Type~$2$ and Type~$1$ $C_4$'s.
For convenience in the following we will bound the number of copies of $C_4$ incident to a given vertex by the sum of five terms where each term represents an upper bound on the number of $C_4$'s of Type~$4$, Type~$3$, red Type~$2$, blue Type~$2$ and Type~$1$ in this given order. 

If $s=2$, then for each pair of vertices there is at most one $C_4$ containing this pair as opposite vertices. 
Therefore we have $C_4(G)\leq \frac{\binom{n}{2}}{2}< \binom{n-2}{2}+2$ for $n>5$.

If $s=3$ and $G$ contains a $K_5$, then
without loss of generality we may assume $G[V(K)]$ is isomorphic to $K_5$.
Then every $C_4$ incident to a vertex of $K$ is either of Type~$4$ or Type~$1$ since $G$ is $C_6$-free.
Note that every vertex of $G-K$ is adjacent to at most one vertex of $K$, hence by the pigeonhole principle there is a vertex $v$ of $K$ with at most $\frac{n-5}{5}$ neighbors in $V(G-K)$.
Hence by applying Lemma~\ref{lemma:c_4vertex} for the vertex $v$ in $G[V(G-K)\cup \{v\}]$, we have that the number of copies of $C_4$ incident to $v$ in $G$ is at most
\[
12+0+0+0+3\cdot\frac{n-5}{5}< n-3,
\]
and the last inequality holds for $n>30$.

If $s=3$ and $G$ is $K_5$-free,
then there are at least three vertices of $K$ not contained in any blue Type~$2$ $C_4$ by Observation~\ref{Obs:general_observation_red_blue}.
Hence by the pigeonhole principle, there is a vertex $v$ of $K$ not contained in a blue Type~$2$ $C_4$ with at most $\frac{n-5}{3}$ neighbors in $V(G-K)$. 
Thus by applying Lemma~\ref{lemma:c_4vertex} and Observation~\ref{Obs:general_observation_red_blue} for the vertex $v$ in $G[V(G-K)\cup \{v\}]$, we have that the number of copies of $C_4$ incident to $v$ is at most
\[
12+2+1+0+2\cdot \frac{n-5}{5}< n-3,
\]
and the last inequality holds for $n\geq 27$.

If $s\geq 4$ and there is a red edge, then there is no blue edge incident to any of the vertices of $B$.
Let us fix a vertex $v$ in $B$ not incident to the red edge and with a neighborhood of the smallest size in $V(G-K)$.
Then by the pigeonhole principle, we have that $\abs{N_{G}(v)\cap V(G-K) } \leq\frac{n-s-2-1}{s-1}$ since there are $n-s-2$ vertices in $G-K$ and at least one is incident to the red edge.
Since there are no red or blue edges incident to the vertex $v$, there are no Type~$2$ and Type~$3$ $C_4$'s incident to $v$. 
Hence we have that the number of copies of $C_4$ incident to $v$ is at most
\[
(s-1)+0+0+0+  \max \left(3\cdot \frac{n-s-3}{s-1}, \frac{(s-1)(s+2)}{2(s+1)}\cdot \frac{n-s-3}{s-1}   \right)<n-3
\]
for all $n>6$.

Assume $4\leq s<\frac{n}{3}$ and there is no red edge.
Let us fix a vertex $v$ in $B$ with a neighborhood of the smallest size in $V(G-K)$. 
Then by the pigeonhole principle, we have~$\abs{N_{G}(v)\cap V(G-K) } \leq\frac{n-s-2}{s}$. 
Thus by applying Lemma~\ref{lemma:c_4vertex} and Observation~\ref{Obs:general_observation_red_blue} for the vertex $v$ in $G[V(G-K)\cup \{v\}]$, we have that the number of copies of $C_4$ incident to $v$ is at most
\[
(s-1)+0+0+(s-1)+\max \left(3\cdot \frac{n-s-2}{s},  \frac{(s-1)(s+2)}{2(s+1)}\cdot \frac{n-s-2}{s}   \right)<n-3
\]
for all $n>18$.

If $\frac{n}{3}\leq s<n-2$ and there is no red edge in $K$, then there is a vertex $v$ in $B$ not incident to a blue $C_4$ such that it has at most  one neighbor   in $V(G-K)$ . 
Therefore the number of $C_4's$ incident to $v$ is at most $(s-1)+0+0+0+0<n-3$.

If $s=n-2$, then $G$ contains a copy of $K_{2,n-2}$. In this case we have $C_4(G)\leq \binom{n-2}{2}+2$, and equality holds only for a graph obtained from a graph isomorphic to $K_{2,n-2}$ by adding an edge in each independent set.

For all $n>3\binom{31}{4}$, we have either $C_4(G)\leq \binom{n-2}{2}+2$ and equality holds  only for a graph obtained from  $K_{2,n-2}$  by adding an edge in each independent set, or
\[
C_4(G)\leq 3\binom{31}{4}+\sum_{i=31}^{n}(i-4)<\binom{n-2}{2}.\qedhere
\]
\end{proof}

\subsection{Maximizing complete bipartite graphs}

\begin{proof}[Proof of Theorem~\ref{thm:kss}]
Since $K_{s,n-s}$ contains $\binom{n-s}{s}$ copies of $K_{s,s}$, we have~$\ex(n,K_{s,s},C_{2s+2})\geq \binom{n-s}{s}$.
In the following we prove a matching upper bound for all integers $n$ and $s$ such that $s\geq 4$ and $n\geq 3s+2+\binom{2s+1}{s}\frac{s+1}{2}$. 
Let $G$ be an $n$-vertex $C_{2s+2}$-free graph with the maximum number of copies of $K_{s,s}$. 

Let $K$ be a subgraph of $G$ isomorphic to $K_{s,s}$ with color classes $A$ and $B$. 
If there is a subgraph $K'$ of $G$ isomorphic to $K_{s,s}$ with color classes $A'$ and $B'$ different from $K$ such that $A=A'$, then there exists a vertex $b\in B'\setminus B$ incident to all vertices of $A$.
Similarly if there is a subgraph $K''$ of $G$ isomorphic to $K_{s,s}$ with color classes $A''$ and $B''$ different from $K$ such that $B=B''$, then there exists a vertex $a\in A''\setminus A$ incident to all vertices of $B$. 
Since $G$ is $C_{2s+2}$-free, we have $a=b$ and $K$ is a contained in a  subgraph of $G$ isomorphic to $K_{s,s}+K_1$. 
Thus for each subgraph $K$ of $G$ isomorphic to $K_{s,s}$, either there is a color class $A$ of $K$ such that the only $K_{s,s}$ containing $A$ as a color class is $K$ or there is a vertex in $G$ incident to all vertices of $K$.
Observe that the graph $K_{s,s}+K_1$ contains cycles of every length $l\in \{3,4,\dots,2s+1\}$.

In the following we give a classification of all subgraphs of~$G$ isomorphic to~$K_{s,s}$.

\begin{itemize}
    \item  For each subgraph $M$ of $G$ such that $M \cong K_{s,s}+K_1$, let us label all copies of $K_{s,s}$ which are subgraphs of $G[V(M)]$ by $G[V(M)]$. 
    We call each such label a \emph{Type 1} label.
    
    \item For each maximal subgraph $M$ of $G$ such that $M \cong K_{s,t}$ for some $t\geq s+2$ let us label all copies of $K_{s,s}$ which are subgraphs of $G[V(M)]$ by $G[V(M)]$.
    We call each such label a \emph{Type 2} label.
    
    \item Let $M$ be a maximal subgraph of $G$ such that $M \cong K_{s,t}$ for some $t=s$ or $t=s+1$ and $K_{s,s}+K_1 \not\subseteq G[V(M)]$. 
    If $t=s+1$, let us label all copies of $K_{s,s}$ which are subgraphs of $G[V(M)]$ by $G[V(M)]$. If $t=s$, let $V'$ be the set of vertices from  $V(G-M)$ incident to a vertex from both color classes of $M$ and incident to at least two vertices from some color class of~$M$. 
    Note that $\abs{V'}\leq 1$ since $G$ is $C_{2s+2}$-free.
    We label each copy of $K_{s,s}$ which is a subgraph of $G[V(M)\cup V']$ by $G[V(M)\cup V']$, and we call each such label a
    \emph{Type 3} label.
\end{itemize}
 
It is easy to see that each copy of $K_{s,s}$ is assigned at least one label. 
Even more, we will now show that each copy of $K_{s,s}$ in $G$ receives exactly one label. 
If a copy $K$ of $K_{s,s}$ has a Type~1 label $G[V(M)]$, then it is immediate from the definition of the labels that it has no Type~3 label. 
Since $G$ is $C_{2s+2}$-free, there is no vertex from $V(G-M)$ incident to at least two vertices of $K$.
Therefore $K$ has no Type~1 label distinct from $G[V(M)]$ and no Type~2 label either. 
By the definition of Type~2 and Type~3 labels, it is impossible to have both types of labels for a copy of $K_{s,s}$. 
If a copy of $K_{s,s}$ received a Type~2 label, then it contains a color class $A$ such that the only $K_{s,s}$ containing $A$ as a color class is this $K_{s,s}$, hence no $K_{s,s}$ receives two different Type~2 labels. 
Finally, no copy of $K_{s,s}$ receives two distinct Type~3 labels, say $G[V(K)\cup\{v_1\}]$ and $G[V(K)\cup\{v_2\}]$, for otherwise  $G[V(K)\cup\{v_1,v_2\}]$ contains a cycle of length $2s+2$, a contradiction.

\begin{claim}\label{Claim_two_Kss_color_class_claim}
Let $K$  and $K'$ be two copies of $K_{s,s}$ in $G$ with color classes $A,B$ and $A',B'$, respectively, such that $K$ and $K'$ have different labels. 
Then we have $A\not\subseteq A' \cup B'$ and $B \not\subseteq A' \cup B'$.
\end{claim}
\begin{proof}
Suppose by way of contradiction that $A\subseteq A' \cup B'$. First note that since the labels are distinct, neither $A=A'$ nor $A=B'$. 
Hence we may assume without loss of generality that $|A\cap A'|\geq 2$ and $|A \cap B'|\geq 1$.
Thus there are at most $2s+1$ vertices in $V(K)\cup V(K')$, otherwise there would be at least two vertices in $B\setminus (A'\cup B')$ incident to two vertices in the class $A'$ and to a vertex in $B'$, yielding a cycle of length $2s+2$.
Note that $V(K)\neq V(K')$ since they have different labels. 
Hence we have $\abs {V(K)\cup V(K')}=2s+1$.
We denote the only vertex in $V(K)\setminus V(K')$ by $v'$. 
This vertex has at least two neighbors in the color class $A'$, which we denote by $a_1, a_2$, and at least one neighbor in the color class $B'$, which we denote by $u$.

Since $K'$ has a different label than $K$ and $v'$ has at least two neighbors in $A$ and a neighbor in $B$, there is a vertex~$v''$ from $V((G-K')-K)$ incident to all vertices of~$A'$ or all vertices of~$B'$.
 Since $G$ is $C_{2s+2}$-free and $v'$ has at least two neighbors in $A'$, the vertex $v''$ is incident to all vertices of $A'$.
 Even more, $v'$ has exactly one neighbor in $B'$, which was previously denoted by $u$. 
 Observe that $B\cap B'$ is nonempty; let us denote an arbitrary vertex from $B\cap B'$ by $u'$. 
 Since $A'$ and ($B'\setminus \{u\})\cup \{v''\}$ induce a complete bipartite graph, there exists a cycle of length $2s$ using the edge  $u'a_1$; let us denote this cycle by $C$. 
 Note that there is a path $P=u'uv'a_1$  of length three.
 By replacing the edge $u'a_1$ of $C$ with $P$, we obtain a cycle of length $2s+2$, a contradiction.
 \end{proof}

\begin{claim}\label{claim:2_vertex_inbtersection}
Let $K$  and $K'$ be two copies of $K_{s,s}$ in $G$ with color classes $A,B$ and $A',B'$, respectively, such that $K$ and $K'$ have different labels. Then we have either
\[
\abs{V(K)\cap V(K')}\leq 1,
\]
or $\abs{V(K)\cap V(K')}=2$ and the vertices of $V(K)\cap V(K')$ are in the same color class in one of the copies of $K_{s,s}$ and in different color classes in the other copy of $K_{s,s}$.
\end{claim}

\begin{proof}
First we prove $\abs{A\cap A'}\leq 1$.
Suppose by way of contradiction that $\{a,a'\}\subseteq A\cap A'$, where $a\neq a'$.
Then by Claim~\ref{Claim_two_Kss_color_class_claim}, there exist vertices $a_1$ and $b_1$  such that $a_{1}\in A'\setminus V(K)$ and $b_1 \in B'\setminus V(K)$. 
Since $B'$ does not contain the vertices $a$ and $a'$,  there is a vertex $b_2 \in B'\setminus (A\cup\{b_1\})$. 

Note that there exists a $2s$-cycle $C$ in $K$ such that $a$ and $a'$ are at distance two. Even more, if $b_2\in B$, then let $C$ be a $2s$-cycle in which the vertex $b_2$ is incident to $a$ and $a'$. Then by changing the path of length two in $C$ between the vertices $a$ and $a'$ with the path $a b_1 a_1 b_2 a'$, we obtain a cycle of length $2s+2$, a contradiction. 

In a similar way, we have that if $\abs{A\cap A'}=1$, then  $\abs{B\cap B'}=0$. This implies the claim.
\end{proof}

Let us fix a maximal copy $M$ of $K_{s,t}$ in $G$ with color classes $D$ of size $s$ and $D'$ of size $t$. 
If $s=3$, we may assume that $t>4$, otherwise for every pair of vertices $u,v$ there are at most~$20$ copies of $K_{3,3}$ containing $v$ and $u$ in the same color class by Claim~\ref{claim:2_vertex_inbtersection}, since the number of $K_{3,3}$'s in a $K_7$ containing a fixed pair of vertices in the same color class is $20$. Then 
\[
K_{3,3}(G)\leq \frac{20}{6}\binom{n}{2}<\binom{n-3}{3}
\]
since $n$ is sufficiently large.

In what follows, we find an injective mapping $\phi$ from the set of copies of $K_{s,s}$ in $G$ to $\binom{V(G)\setminus D}{s}$. Finding such a mapping implies the desired upper bound. 

Let $L_1$ be a copy of $K_{s,s}+K_1$ in $G$, if such a copy exists. Note that $G[V(L_1)]$ is a Type~1 label. If $s>3$ and  $V(L_1)\neq V(M)$, then let $A'$ be a subset of $V(L_1)\setminus M$ of size~$s-1$. Such a subset exists by Claim~\ref{claim:2_vertex_inbtersection}. If $s>3$ and $V(L_1)= V(M)$, let $A'$ be a subset of $V(L_1)\setminus D$ of size~$s-1$.
Let~$X$ be a set of all copies of $K_{s,s}$ in $G[V(L_1)]$. Note that $\abs{X}\leq K_{s,s}(K_{2s+1})= \frac{s+1}{2}\binom{2s+1}{s}$. Let $Y$ be a subset of $V(G)\setminus (V(L_1)\cup D)$ of size $\abs{X}$; such a subset exists since $n$ is large. We define $\phi$ on the set of copies of $K_{s,s}$ with label $G[V(L_1)]$ in the following way
\[
\phi(X):=\left\{ A'\cup \{v\}:v\in Y\right\}. 
\]
If $s=3$, then  $V(L_1)\neq V(M)$ since $t>4$. Then by Claim~\ref{claim:2_vertex_inbtersection}, we have $\abs{V(L_1) \cap D}\leq 1$   and $\abs{V(L_1) \cap D'}\leq 1$. Let $X$ be a set of all copies of $K_{s,s}$ in $G[V(L_1)]$. 
Note that $|X|\leq 70$. 
If $\abs{V(L_1) \cap D}=0$, then let $Y$ be a set of size $\abs{X}$ containing distinct triples from  $\binom{V(L_1)}{3}$ and  $\binom{V(L_1)}{2}$ with a vertex from $D'\setminus V(L_1)$. Note that $\abs{D'\setminus V(L_1)}\geq 2$ by Claim~\ref{claim:2_vertex_inbtersection}. Such a set~$Y$ exists since $|D' \setminus V(L_1)| \ge 2$ and $\abs{X}\leq 70<\binom{7}{3}+2\binom{7}{2}$. 
If $\abs{V(L_1) \cap D}=1$, then let $Y$ be a set of size $\abs{X}$ containing distinct triples from $\binom{V(L_1)\setminus D}{3}$ and  $\binom{V(L_1)\setminus D}{2}$ with  vertices  from $D'\setminus V(L_1)$. 
Note that $\abs{D'\setminus V(L_1)}\geq 4$ by Claim~\ref{claim:2_vertex_inbtersection} and $70<\binom{6}{3}+4\binom{6}{2}$.
We define $\phi$ on the set of copies of $K_{s,s}$ with label $G[V(L_1)]$ by
\[
\phi(X):=\left\{ A'\cup \{v\}:v\in Y\right\}. 
\]

Let $L_2$ be a maximal subgraph of $G$ isomorphic to $K_{s,t'}$ for some $t'\geq s+2$ with color classes $A$ of size $s$ and $B$ of size $t'$, if such a subgraph $L_2$ exists. Note that $G[V(L_2)]$ is a Type 2 label.
Since $G$ is $C_{2s+2}$-free, each  subgraph $K$ of $G[V(L_2)]$ isomorphic to $K_{s,s}$ has color classes $A$ and  a subset $B'$ of $B$ of size $s$. 
Then by Claim~\ref{claim:2_vertex_inbtersection}, we have $\abs{B \cap D}\leq 1$.
If $\abs{B \cap D}=0$, then let $\phi(K):=B'$ and if $\abs{B \cap D}\leq 1$, then let $\phi(K):=(B\cup\{v_{L_2}\})\setminus D$, where $v_{L_2}$ is any fixed vertex of $D' \setminus V(L_2)$.

Let $L_3$ be a subgraph of $G$ such that $G[V(L_3)]$ is a Type~3 label if such a subgraph $L_3$ exists.  We define $\phi$ for the copies of~$K_{s,s}$ in~$G[V(L_3)]$ as we did for all copies of~$K_{s,s}$ with Type~1 label. Note that if $s=3$ and $\abs{V(L_3)\setminus D}\leq 1$ the number of $K_{3,3}$ with such a label is at most $10$, and the number of triples $\binom{V(L_3)\setminus D}{3}$ is at least $10$. Therefore for the image of $\phi$ we do not need to use any vertex outside of this label.
If $s=3$ and $\abs{V(L_3)\setminus D}=2$ (by Claim~\ref{claim:2_vertex_inbtersection} it is at most $2$), then 
$L_3$ induces a $K_{3,3}$ since $G$ is $C_8$-free. Therefore for the image of $\phi$ we do not need to use any vertex outside of this label.

Let $K$ and $K'$ be copies of $K_{s,s}$ in $G$.
Note that if $K$ and $K'$ have the same label, then we have 
\[
\phi(K) \neq \phi(K').
\]
If   $K$ and $K'$ have different labels and $s>3$,  then  since each triple $\phi(K)$ and $\phi(K')$ contains at most one vertex outside of its label and we have $\phi(K) \neq \phi(K')$ by Claim~\ref{claim:2_vertex_inbtersection}. 

If $K$ and $K'$ have different labels and $s=3$ suppose by way of contradiction that $\phi(K)=\phi(K')=\{a,d,d'\}$. 
By Claim~\ref{claim:2_vertex_inbtersection} we have that $a$ is a vertex of both labels, $d$ is a vertex of label of $K$,  $d'$ is a vertex of label of $K'$  and $d,d'\in D'$. 
There is a path $dPd'$ of length four going through the vertex $a$ in $G[V(K)\cup V(K')]$.
There is a path of length $5$ in $G[M]$ from $d$ to $d'$ internally disjoint from $dPd'$ by Claim~\ref{claim:2_vertex_inbtersection}. 
These two paths constitute a cycle of length $8$, a contradiction. 
Hence~$\phi$ is an injective function and we obtained the extremal number $ex(n,C_6,C_8)$.

Note that if there is more than one label, then for each label there is room to choose an extra set from $\binom{V(G)\setminus D}{s}$ which will not be an image of $\phi$ at the end of the procedure. 
Therefore equality holds only for graphs~$G$ containing exactly one label. 
Hence $K_{s,n-s}\subseteq G$. 
It is easy to see that $G\subseteq K_s+H$, where~$H$ is an $(n-3)$-vertex graph with exactly one edge since $G$-is $C_8$-free.
\end{proof}

\section{Acknowledgements}
 The research of Gy\H{o}ri and Salia was supported by the National Research, Development and Innovation Office NKFIH, grants  K132696 and SNN-135643. The research of Tompkins was supported by NKFIH grant K135800.

\bibliography{references.bib}

  \textit{E-mail addresses:} \\
  E.~Gy\H{o}ri: \texttt{gyori.ervin@renyi.hu}\\
  Z.~He: \texttt{hz18@mails.tsinghua.edu.cn}\\
  J.~Lv: \texttt{lvzq19@mails.tsinghua.edu.cn}\\
  N.~Salia: \texttt{nikasalia@yahoo.com}\\
  C.~Tompkins: \texttt{ctompkins496@gmail.com}\\
  K.~Varga: \texttt{vkitti@renyi.hu}\\
  X.~Zhu: \texttt{ zhuxt@smail.nju.edu.cn}

\end{document}